\documentclass[journal,twoside]{IEEEtran}

\usepackage{amssymb,amsmath,color,graphicx, amsbsy}
\usepackage{graphicx}
\usepackage{subfigure}
\usepackage{cite}
\usepackage{algorithm}
\usepackage{algorithmic}
\usepackage{amsfonts}
\usepackage{amssymb}
\usepackage{graphicx}
\usepackage{epstopdf}
\usepackage{epstopdf}

\newtheorem{theorem}{Theorem}%[section]

\newtheorem{definition}[theorem]{Definition}

\ifodd 1
\newcommand{\rev}[1]{{#1}}
\newcommand{\com}[1]{{\color{red} #1}}
\newcommand{\res}[1]{\textbf{\color{magenta} (RESPONSE: #1)}}
\newcommand{\del}[1]{}
\else
\newcommand{\rev}[1]{#1}
\newcommand{\com}[1]{}
\newcommand{\res}[1]{}
\newcommand{\del}[1]{}
\fi
\begin{document}
\title{\rev{Risk-limiting} Economic Dispatch for Electricity Markets with Flexible Ramping Products}

\author{Chenye~Wu,~\IEEEmembership{Member,~IEEE,}
        Gabriela Hug,~\IEEEmembership{Senior Member,~IEEE,}
        and~Soummya Kar,~\IEEEmembership{Member,~IEEE}\\ \vspace{-0.5cm}
\thanks{C.\,Wu, G. Hug and S. Kar are with the Department of Electrical and Computer Engineering, Carnegie Mellon University, Pittsburgh, PA 15213, USA, e-mails: chenyewu@andrew.cmu.edu, ghug@ece.cmu.edu, soummyak@andrew.cmu.edu.
} }

\maketitle
\thispagestyle{empty}
\pagestyle{empty}

\begin{abstract}
The expected increase in the penetration of renewables in the approaching decade urges the electricity market to introduce new products - in particular, \emph{flexible ramping products} - to accommodate the renewables' variability and intermittency. \rev{A risk-limiting economic dispatch scheme provides the means to optimize the dispatch and provision of these products. In this paper, we adopt the extended loss-of-load probability as the definition of risk. We first assess how the new products distort the optimal economic dispatch by comparing to the case without such products. Specifically, using parametric analysis, we establish the relationship between the minimal generation cost and the two key parameters of the new products: the up- and down-flexible ramping requirements.} Such relationship yields a novel routine to efficiently \rev{solve the non-convex risk-limiting economic dispatch problem}. Both theoretical analysis and simulation results suggest that our approach may substantially reduce the cost for incorporating the new products. We believe our approach can assist the ISOs with utilizing the ramping capacities in the system at the minimal cost.
\end{abstract}

%\vspace{-0.3cm}

\begin{keywords}
Flexible ramping product, real time market, renewable energy integration, parametric optimization
\end{keywords}

%----------------------------------------------------------------------
% SECTION I: Introduction
%----------------------------------------------------------------------
\section{Introduction}
\label{intro}

As the level of penetration of renewable generation (in particular wind and solar power) grows, the stochastic nature of the power outputs from these resources is increasingly stressing the power system. Hence, the NERC task force on the potential reliability impacts of emerging flexible resources \cite{NERCtaskforce} suggests designing new products for the future electricity market. These products should ensure sufficient ramping capacity to cope with large short-term variations and prediction errors \cite{lannoye2010integration}. Activity is starting to pick up: the California Independent System Operator (CAISO) \cite{abdul2012enhanced} and the Midwest Independent System Operator (MISO) \cite{navid2012market} are pioneering the design and have introduced flexible ramping products in their markets.

Conceptually, flexible ramping products aim at reserving ramping flexibility in the current time slot for future use. While frequency regulation already reserves certain flexibility to tackle unpredicted fluctuations in net load \cite{bergenpower}, the new products are expected to provide more flexibility than frequency regulation and on a much slower time scale (e.g., 5 minutes in CAISO). Spinning and non-spinning reserves are other existing products to provide flexibility but are held to manage system contingencies. Furthermore, they can only contribute to up-ramping flexibility while flexible ramping products are intended to tackle deviations in both directions.

%Conceptually, these flexible ramping products aim at providing more ramping capacity (flexibility) to the system by reserving capacity in the current time slot for future use. This feature contrasts it from frequency regulation, spinning reserve, and non-spinning reserve, since the existing products are held to manage certain pre-defined system contingencies occurring in the same interval. Also, note that, spinning reserve and non-spinning reserve can only contribute to up-ramping flexibility while the flexible ramping products should cope with the deviation in both directions. Furthermore, different from frequency regulation, the flexible ramping products are designed to achieve a cost effective dispatch in 5-min intervals whereas frequency regulation cares less about the cost effectiveness but puts an emphasis on the capability to react to the control signal (automatic generation control signal) within several seconds \cite{bergenpower}. The amount of flexibility provided by the new products is also expected to be much higher than that of frequency regulation.

Albeit promising and important, \rev{to better utilize these products, a risk-limiting economic dispatch scheme is warranted, in which we adopt the extended loss-of-load probability (LOLP) \cite{billinton1970power} as the definition of risk. The risk-limiting economic dispatch scheme is in general non-convex due to the constraints to limit the risk. This non-convexity heavily constrains the risk-limiting economic dispatch scheme from (near) real time implementation. We propose an efficient algorithm, which utilizes our understanding on} how the new products distort the optimal economic dispatch compared to the market outcomes for the case without such products. The understanding is motivated by the following question: \rev{what are the main factors that determine the distortion}?

Intuitively, it depends on how much ramping flexibility is needed. This value, on the other hand, is dependent on the desired level of reliability at which the load in the next time step should be supplied. \rev{In this paper, we adopt the risk-limiting constraint to enforce the desired level of supply reliability (or equivalently, LOLP).} As different combinations of the up- and down-flexible ramping capacities may meet the same risk-limiting constraint, it is possible to optimize the combination by minimizing the total generation cost. \rev{ We propose to solve this problem in two steps. First,} %In general, however, such kind of optimization is extremely difficult in that it requires a brute-force search and needs significant efforts on solving large amount of optimization problems.
by employing the theory of linear parametric programming \cite{Holder2010}, we propose a parametric functional approach to understand the dependency of the distortion cost on the level of up- and down-ramping requirements. Then, we introduce a linear search algorithm to solve the economic dispatch while guaranteeing a pre-defined risk-limiting constraint. %Such an approach grants us the possibility to carry out the brute-force search efficiently.

%On the other hand, the product has two key parameters: the up and down flexible ramp requirements. If we are able to select different parameters to achieve the same level of system reliability, can we achieve a more effective market by reducing the distortion cost? With the increasing penetration level, the two flexible ramp requirements will increase correspondingly. Hence, a natural research question is to quantify the relationship between the distortion cost and the increasing flexible ramp requirement, and finally the increasing penetration level of renewable generation.
%%
%To address all of these questions, in this paper, we propose a functional approach to analytically understand the relationship between the distortion cost and key parameters of the flexible ramping products. By considering the physical laws and the power network constraints, we also make several interesting observations on the role of power line congestion in achieving a cost effective market.

It is worth noting that the application of the proposed parametric analysis is not limited to the cost assessment. Examples of such application include: the parametric optimal power flow (OPF) to offer an excellent visualization of the complex nature of OPF \cite{almeida1994general}, and recently the unifying functional approach to assessing the market power \cite{boseunifying}.

%\vspace{-0.2cm}

\section{Related Work and Contributions}

\subsection{Related Work}
The research on flexible ramping products started only recently. Wang \emph{et al.} performed extensive comparisons between the market outcomes of the ISO's deterministic market model and the optimal stochastic model for flexible ramping products in \cite{wang2014flexible}. Taylor \emph{et al.} designed an optimal dynamic pricing scheme for the ancillary service markets including the flexible ramping products in \cite{taylor2013dynamic}. Our paper further this track by understanding the relationship between the flexible ramping capacity requirements and the electricity market outcomes.

\rev{Our work also fits into a growing literature on the risk-limiting economic dispatch. For example, Varaiya \emph{at al.} introduced a conceptual framework for risk-limiting economic dispatch in \cite{5618534}. Rajagopal \emph{et al.} furthered the research by proposing a closed-form computational model in \cite{Rajagopal2013615}. Zhang \emph{et al.} utilized the Monte Carlo sampling based scenario approximation technique to conduct the risk-limiting economic dispatch in \cite{6497884}. Different from the previous work, we highlight the impact of flexible ramping products in risk-limiting economic dispatch and propose an efficient algorithm to solve the problem.}

Another set of related work focused on the analysis of the cost brought by the variability of renewables. For example, Katzenstein \emph{et al.} introduced a novel metric for evaluating the cost of wind power variability in \cite{katzenstein2012cost}. Lueken \emph{et al.} presented the costs induced by the solar and wind power in \cite{lueken2012costs}. In contrast to \cite{katzenstein2012cost,lueken2012costs}, we introduce a functional approach to assessing the cost brought by the variability of renewables, and we focus on the distortion cost.

In our earlier work \cite{wuinsubmission}, we made the first step towards understanding flexible ramping products' influence on the electricity market outcomes. In this paper, we further the research by introducing a novel routine to efficiently construct the parametric functions for computational purpose. \rev{Based on this routine, we show how to solve the risk-limiting economic dispatch efficiently with a linear search.}

\vspace{-0.1cm}
\subsection{Our Contributions}
Towards understanding \emph{the relationship between the generation cost and the key parameters of flexible ramping products, and based on this dependency, how to achieve the minimal generation cost with limited risk}, the major contributions of this paper are summarized as follows:

\vspace{0.05cm}
\begin{itemize} \itemsep 2pt
  \item \emph{Parametric Analysis:} We employ a parametric functional approach to studying the relationship between the generation cost and the up- and down-ramping requirements. Such an approach displays promising properties (such as monotonicity, convexity/concavity, and piecewise linearity) of the function.
  \item \emph{Triple Optimality Guarantee:} Inspired by the parametric analysis, we consider the cost minimization problem from two additional viewpoints: given a certain financial budget and the required up (down) flexible ramping capacity, what is the maximal down (up) flexible ramping capacity that the system can provide? We prove that certain inverse function relationships exist among the three proposed functions, and each of them enjoys triple optimality.
  \item \emph{3D Function Efficient Construction:} \rev{Each of the parametric optimization functions (formally defined in Section \ref{functional}) has two arguments, which is in general hard to construct efficiently. By utilizing the triple optimality guarantee, we propose two efficient routines to construct the parametric functions: one for computation, the other for visualization.}
  \item \emph{Risk-limiting Cost Minimization:} Based on the proposed efficient function construction, we carry out the linear search for the up- and down-ramping requirements which minimize the generation cost while ensuring the pre-defined risk-limiting constraint (i.e., supply reliability constraint). \rev{This essentially addresses the non-convex risk-limiting economic dispatch problem.}
\end{itemize}
\vspace{0.05cm}

The rest of this paper is organized as follows: we revisit the mathematical formulation for incorporating these products into the real time energy dispatch market, and \rev{highlight the notion of risk-limiting economic dispatch} in Section \ref{system}. \rev{After identifying the challenges to solve the risk-limiting economic dispatch problem, we first closely investigate the classical economic dispatch problem (without the risk-limiting constraint). }
%we try to understand the optimization problem from two different viewpoints:
%\begin{itemize}
%  \item the impact of the parameters on the distortion cost;
%  \item given a certain financial budget and the required up (down) ramping capacity, the maximal down (up) flexible ramping that the system can provide.
%\end{itemize}
By applying the parametric functional approach from different aspects to the classical economic dispatch problem, we propose three parametric optimization functions in Section \ref{functional}. Subsequently, in Section \ref{analytical}, we investigate various analytical properties of these functions to draw a clear relationship between the generation cost and the ramping requirement parameters. Based on these relationships, we \rev{revisit the risk-limiting economic dispatch problem and propose a linear search algorithm to determine the two key parameters for the flexible ramping products} in Section \ref{distortion}. Section \ref{simulation} presents several illustrative examples and case studies to evaluate the performance of our approach. Finally, our concluding remarks and directions for future work are discussed in Section~\ref{conclusions}.

%----------------------------------------------------------------------
% SECTION II: System Model
%----------------------------------------------------------------------
%\vspace{-0.2cm}
\section{Problem Formulation}
\label{system}

CAISO proposes implementing flexible ramping products in the 5-minute real time market. Hence, in this paper, we first cast the problem as a \emph{model predictive control} (MPC) problem \cite{camacho2013model} with a horizon of $T$ time steps, each being of 5 minute length. \rev{This formulation can be easily generalized to other markets by selecting proper length of the time scale.} The objective is to minimize generation cost to supply the expected load subject to the DC load flow constraints, the limitations on generation outputs, and \rev{the limitations on risk with respect to supply reliability by determining the required levels of generation ramping capacity. For notational simplicity, we do not consider other products, e.g., frequency regulation, spinning reserves, and non-spinning reserves, in this model. In fact, they will only incur linear constraints, which will not affect our subsequent analytical results. Mathematically, the risk-limiting economic dispatch problem with flexible ramping products can be formulated as follows:}
%In our model, we consider the  real time economic dispatch including flexible ramping requirements. The set of all generators is denoted by $\mathcal{N}$. %Suppose each generator $n\in\mathcal{N}$ submits a linear bid $C_n$ into the energy market, then
%The ISO may seek to minimize the generation cost by performing the following optimization:
\begin{align}
   \min \ &  \sum_{t=0}^{T-1} \sum_{n\in\mathcal{N}} \left ( C_n g_{n,t}\! + \rev{k_n^u r^u_{n,t} + k_n^d r^d_{n,t}} \right ) \label{op_obj}\\
   \text{{\emph{s.t.} }} & -\boldsymbol{b} \le H_g \boldsymbol{g}_t - H_d\boldsymbol{\hat{d}}_t \le \boldsymbol{b}, \ \forall t, \label{op_linecapacity} \\
   & \boldsymbol{1}^T \boldsymbol{g}_t - \boldsymbol{1}^T \boldsymbol{\hat{d}}_t = 0, \ \forall t, \label{op_total} \\
   %& \underline{g}_n \le g_{n,t} \le \bar{g}_n, \ \forall n, \ \forall t, \label{op_capacity} \\
   & \textstyle \sum_{n\in\mathcal{N}} r^u_{n,t} = F^{u}_t, \  t=1,\cdots,T-1, \label{op_rampup_capacity} \\
   & \textstyle \sum_{n\in\mathcal{N}} r^d_{n,t} = F^{d}_t, \  t=1,\cdots,T-1, \label{op_rampdown_capacity} \\
  % & -\Delta g_n \le -r^d_{n,t} \le 0, \ \forall n, \ \forall t, \label{op_rampdown} \\
   & \underline{g}_n \le g_{n,t} + r^u_{n,t} \le \bar{g}_n, \ \forall n, \ \forall t, \label{op_capacity_up} \\
   & \underline{g}_n \le g_{n,t} - r^d_{n,t} \le \bar{g}_n, \ \forall n, \ \forall t, \label{op_capacity_down} \\
   % & 0 \le r^u_{n,t}, \ r^d_{n,t} \le \Delta g_n, \ \forall n, \ \forall t, \label{op_rampup} \\
   &  |g_{n,t+1} - g_{n,t}+ r^u_{n,t+1} + r^d_{n,t}| \le \Delta g_n, \ \forall n,  \forall t, \label{op_maxrampup} \\
   &  |g_{n,t+1} - g_{n,t} - r^d_{n,t+1} - r^u_{n,t}| \le \Delta g_n, \ \forall n,  \forall t, \label{op_maxrampdown}\\
   &  g_{n,t}\ge 0,  \ \forall n,  \forall t, \label{op_nonnegative}\\
   &  r^u_{n,t}\ge 0, r^d_{n,t} \ge 0, \ \forall n,  t=1,\cdots,T-1, \label{op_nonnegative1}\\
 %  &  F^u_{t}\ge 0, F^d_{t} \ge 0, \ t=1,\cdots,T-1, \label{op_nonnegative2}\\
   & \rev{\Pr ( \boldsymbol{1}^T \boldsymbol{\hat{d}}_t \!-\! F_t^d \! \le \!\boldsymbol{1}^T \boldsymbol{d}_t\! \le\! \boldsymbol{1}^T \boldsymbol{\hat{d}}_t\!+\!F_t^u ) \ge p\%, \forall t.} \label{op_risk}
\end{align}
The decision variables in the problem (\ref{op_obj})-(\ref{op_risk}) are
\begin{itemize}
  \item $g_{n,t}$: generator $n$'s power output [MW] at time $t$, with vector form $\boldsymbol{g}_t = [g_{n,t}, \forall n\in\mathcal{N}]$;
  \item $r^u_{n,t}$, $r^d_{n,t}$: up and down flexible ramping capacities [MW] provided by generator $n$ at  time $t=1,\cdots,T-1$;
  \item $F^{u}_t$, $F^{d}_t$: up- and down-ramping requirements [MW] for the overall system at time $t=1,\cdots,T-1$;
\end{itemize}
%Note that $g_{n,0}$'s are considered as given values, \com{and
%$r^{u}_{n,0}=r^{d}_{n,0}=0,\ \forall n\in\mathcal{N}.$ The proposed MPC approach seeks to perform the final economic dispatch for time $t=1$. As shown in Fig. \ref{fig:visual}, though there could possibly be minor fluctuations at time  $t=1$, such fluctuations should be taken care of by frequency regulation, not the flexible ramping products. The ramping capacity for time $t=1$ was reserved by the optimization in the previous time slot. Therefore, we set $r^{u}_{n,1}=r^{d}_{n,1}=0,\ \forall n\in\mathcal{N}$. The MPC approach also needs to reserve certain ramping capacities for time $t=2,\cdots,T$ to tackle the major uncertainties (the maximal and minimal possible load profiles shown in Fig. \ref{fig:visual}). Hence, $r^u_{n,t}$'s and $r^d_{n,t}$'s at time $t=2,\cdots,T$ are the variables.}
and the parameters are

\vspace{0.1cm}
\begin{itemize}
  \item $C_n$: bid [\$/MW] of generator $n$ to provide energy;
  \item \rev{$k_n^u, k_n^d$: bids [\$/MW] of generator $n$ for providing ramping up and ramping down capacities;}
  \item $\boldsymbol{b}$: transmission line capacity vector [MW];
  \item $H_g$, $H_d$: generation and load shift factor matrices;
  \item $\boldsymbol{\hat{d}}_t$: predicted demand vector [MW] at time $t$;
  \item $\boldsymbol{d}_t$: actual demand vector [MW] at time $t$;
  \item $\boldsymbol{1}$: unit column vector of appropriate dimension;
  \item $\underline{g}_n$, $\bar{g}_n$: generator $n$'s minimal and maximal generation capacity [MW];
  \item $\Delta g_n$: generator $n$'s ramping limit [MW/5 minutes];
  \item $p$\%: \rev{probability at which the system operator wants to meet the actual demand at all times $t$.}
\end{itemize}
\vspace{0.1cm}

\begin{figure}[t]
\begin{center}
\includegraphics[width=5.5cm]{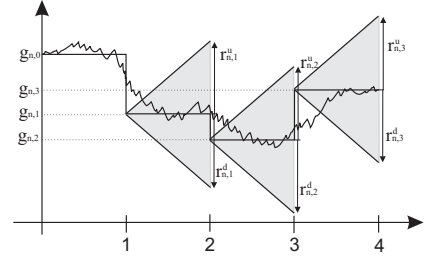}
\vspace{-0.3cm}
 \caption{Visualization of the MPC approach when $T=4$.}\vspace{-0.7cm}
 \label{fig:visual}
\end{center}
\end{figure}

The proposed MPC approach seeks to perform the economic dispatch for time steps $t = 0,...,T-1$ under the condition that ramping capacity needs to be reserved for steps $t=1,...,T-1$. An illustration of the control variables is given in Fig.~\ref{fig:visual}. Ramping capacity for $t=0$ has been reserved in the previous time step, hence, there are no variables $r_{n,0}$ to be determined. Note that the load predictions are updated as time goes by. Hence, only the energy dispatch profile for $t=0$, i.e., $g_{n,0}$'s, and the flexible ramping requirements for $t=1$, i.e., $F_1^u$ and $F_1^d$, will be applied.

Constraint (\ref{op_linecapacity}) corresponds to the line capacity constraints. Constraints (\ref{op_total})-(\ref{op_rampdown_capacity}) represent the total power balance and up and down flexible ramping requirements, respectively. Constraints (\ref{op_capacity_up})-(\ref{op_capacity_down}) ensure that the generation capacity constraints are met and constraints (\ref{op_maxrampup})-(\ref{op_maxrampdown}) enforce that \rev{the ramping limits hold even in the worst cases (i.e., the ramping can be feasibly supplied by the generators if needed)}. The next set of constraints (\ref{op_nonnegative})-(\ref{op_nonnegative1}) ensures that all the decision variables are non-negative. \rev{ The last constraint is the risk-limiting constraint, which implies that the system operator needs to meet the actual demand at all times $t$ with probability of at least $p$\%. Note that in this paper, we regard the renewable energy as negative load. Therefore, all the uncertainties and randomness are in the load predictions, i.e., $\boldsymbol{\hat{d}}_t$'s. Mathematically, we extend the standard probability based risk definition - the Loss-of-Load Probability (LOLP) - found in literature \cite{billinton1970power}:

\vspace{0.1cm}
\begin{definition}
A tuple ($F_t^u,F_t^d$) is said to achieve confidence level of $p$\% with respect to a prediction $\boldsymbol{\hat{d}}_t$ at time $t$, if
\begin{equation}\label{confi_def}
    \Pr ( \boldsymbol{1}^T \boldsymbol{\hat{d}}_t \!-\! F_t^d \! \le \!\boldsymbol{1}^T \boldsymbol{d}_t\! \le\! \boldsymbol{1}^T \boldsymbol{\hat{d}}_t\!+\!F_t^u ) \ge p\%.
\end{equation}
\end{definition}
\vspace{0.1cm}

The probability distribution can be obtained by the prediction error distribution (see Section \ref{errorModelSection} for more details). We want to emphasize that our risk-limiting economic dispatch is different from security-constrained economic dispatch \cite{monticelli1987security}. The latter focuses on network contingent events (e.g., failure of a generator, a transformer, or a line outage). Based on our understanding, in the future, flexible ramping products will be utilized very often. The network contingent events will be handled by spinning reserves, and non-spinning reserves, which are not the focus of this paper.}

In addition, even though we simplify the model in (\ref{op_linecapacity}) and do not capture the influence of flexible ramping products on line flows, they do influence the feasible regions of parameters $F^{u}_t$ and $F^{d}_t$ via the coupling among $g_{n,t}$'s, $r_{n,t}^u$'s, and $r_{n,t}^d$'s in (\ref{op_capacity_up})-(\ref{op_maxrampdown}).

To simplify the subsequent analysis, and to highlight the relationships between the parameters of interest, we concentrate on the analysis for $T=2$. Hence, the only ramping variables are $r_{n,1}^u$ and $r_{n,1}^d$ and we can simplify the notation by setting $F_1^u = F^u$ and  $F_1^d = F^d$.

\rev{Another simplification is to ignore the bids for flexible ramping up and down capacities. This allows us to focus on how the new products distort the optimal economic dispatch. We realize that although in the current implementation, there is no bidding scheme for the new products, in the future, there should be a reasonable bidding scheme to better accommodate the resources. Nevertheless, we want to stress that even without the bidding information, these products will not come free. The prices will be determined by the Lagragian multipliers associated with constraints (\ref{op_capacity_up}) and (\ref{op_capacity_down}). Just as in the case for frequency regulation, this payment is referred to as the capacity payment \cite{federal2011frequency}. One natural question is to examine the relationship between the capacity payment and the generation cost. If there are no line capacity constraints or they are non-binding, by contradiction, we can show that the lower the generation cost is, the lower is the total capacity payment. Although this relationship may not be valid when line capacity constraints become active, the argument possibly still holds over the range of interest. A detailed discussion, however, falls out of the scope of this paper.

Unfortunately, even with all the simplifications, the risk-limiting economic dispatch problem (\ref{op_obj})-(\ref{op_risk}) is still non-convex. Instead of employing a straightforward brute-force search algorithm (which enumerates all the possible solutions) or other iterative algorithms without a guaranteed global optimal solution, we notice that the classical economic dispatch problem (without the risk-limiting constraint) is convex. In particular, with all the simplifications, the resulting economic dispatch problem can be formulated as follows:

\begin{equation}\label{classical_energy_dispatch}
    \begin{aligned}
     \min &  \sum_{t=0}^1  \sum_{n\in\mathcal{N}} C_n g_{n,t} \\
    \text{\emph{s.t.} } &  \textstyle \sum_{n\in\mathcal{N}} r^u_{n,1} = F^{u},\\
    & \textstyle \sum_{n\in\mathcal{N}} r^d_{n,1} = F^{d}, \\
   & \text{Constraints (\ref{op_linecapacity})-(\ref{op_total}), (\ref{op_capacity_up})-(\ref{op_nonnegative1}),}
    \end{aligned}
\end{equation}

Suppose that we are given $F^u$ and $F^d$ as parameters in the simplified economic dispatch problem (\ref{classical_energy_dispatch}), what are the roles of these parameters in determining the minimal generation cost? This question motivates our subsequent parametric functional analysis on the simplified economic dispatch problem. This analysis gives us insights into how to efficiently solve the risk-limiting economic dispatch problem.
}

%Note that the real time economic dispatch is performed every 15 minutes in practice, which renders considering $T>3$ uninteresting. To this end, $T=2$ is already a good practical choice.% for the current electricity market. %Next, we want to investigate the relationship between the distortion cost and these two parameters.

% \vspace{-0.2cm}
\section{Parametric Functional Analysis}
\label{functional}

%With the purpose of investigating how $F^{u}$ and $F^{d}$ influence the distortion cost,
We employ a parametric functional approach to investigating how $F^{u}$ and $F^{d}$ influence the minimal generation cost in the simplified economic dispatch problem (\ref{classical_energy_dispatch}). The key idea is to replace the single-value optimization problem (\ref{classical_energy_dispatch}) with a parameterized function with two parameters. Subsequently, we parameterize another two optimization problems which are closely related to (\ref{classical_energy_dispatch}) and maximize the up and down flexible ramping capacities, respectively, given a certain budget. By fixing certain parameters, we establish the underlying inverse function relationships among the three functions, which imply that each of them enjoys triple optimality - any point on each function in the region of interest is the solution to three optimization problems.

\vspace{-0.2cm}
\subsection{Minimal Cost (MinC) Function}
An extension of the optimization problem (\ref{classical_energy_dispatch}) is to ask \emph{for any given $F^{u}$ and $F^{d}$, what is the minimal generation cost that the ISO could achieve?} This leads to the minimal cost (MinC) function, defined as follows:
\begin{equation}\label{MinC}
    \begin{aligned}
    \text{MinC}(f^{u},f^{d}) = \min &  \sum_{t=0}^1  \sum_{n\in\mathcal{N}} C_n g_{n,t} \\
    \text{\emph{s.t.} } &  \textstyle \sum_{n\in\mathcal{N}} r^u_{n,1} = f^{u},\\
    & \textstyle \sum_{n\in\mathcal{N}} r^d_{n,1} = f^{d}, \\
   & \text{Constraints (\ref{op_linecapacity})-(\ref{op_total}), (\ref{op_capacity_up})-(\ref{op_nonnegative1}),}
    \end{aligned}
\end{equation}
where $f^u$ and $f^d$ are the function's arguments, representing the up and down flexible ramping requirements, respectively. It is immediately clear that MinC$(F^u,F^d)$ corresponds to (\ref{classical_energy_dispatch}). The trivial upper bounds for $f^u$ and $f^d$ are both $\sum_{n\in\mathcal{N}} \Delta g_n$. However, the line capacity constraints and the generation capacity constraints can both shrink the feasible regions of $f^u$ and $f^d$. The lower bounds for $f^u$ and $f^d$ are set to zero, corresponding to the case without any flexible ramping requirement.

%\vspace{-0.2cm}
%\subsection{Distortion Cost (DS) Function}
%\label{distortion_func}
%
%

%\vspace{-0.2cm}
\subsection{Maximal Up Flexible Ramping (MaxUR) Function}
Next, we consider a related optimization problem that the ISO may face. \emph{Given a certain financial budget $\theta$, and the down flexible ramping requirement $f^d$, what is the maximal up flexible ramping that the system can contribute?} Mathematically, we refer to it as the maximal up flexible ramping (MaxUR) function:
\begin{equation}\label{MaxRU}
    \begin{aligned}
    \text{MaxUR}(\theta,f^{d})\! =\!\max & \sum_{n\in\mathcal{N}} r^u_{n,1}  \\
    \text{\emph{s.t.} } &  \!\! \textstyle\sum_{t=0}^1 \sum_{n\in\mathcal{N}}\! C_n g_{n,t}\! \le\! \theta,\\
    & \!\! \textstyle\sum_{n\in\mathcal{N}} r^d_{n,1} = f^{d}, \\
   &\!\! \text{Constraints (\ref{op_linecapacity})-(\ref{op_total}), (\ref{op_capacity_up})-(\ref{op_nonnegative1}).}
    \end{aligned}
\end{equation}
%where $\theta_{b}$ is defined in (\ref{minc00}).

%\vspace{-0.2cm}
\subsection{Maximal Down Flexible Ramping (MaxDR) Function}

Similarly, we can define the maximal down flexible ramping (MaxDR) function by asking, \emph{given a certain financial budget $\theta$, and the up flexible ramping requirement $f^u$, what is the maximal down flexible ramping that the system can contribute?} Mathematically, we have
\begin{equation}\label{MaxDU}
    \begin{aligned}
    \text{MaxDR}(\theta,f^{u}) \!=\!\max & \sum_{n\in\mathcal{N}} r^d_{n,1}   \\
    \text{\emph{s.t.} } & \!\! \textstyle\sum_{t=0}^1 \sum_{n\in\mathcal{N}} \!C_n g_{n,t}\! \le\! \theta,\\
    & \!\! \textstyle\sum_{n\in\mathcal{N}} r^u_{n,1} = f^{u}, \\
   & \!\!\text{Constraints (\ref{op_linecapacity})-(\ref{op_total}), (\ref{op_capacity_up})-(\ref{op_nonnegative1}).}
    \end{aligned}
\end{equation}
%where $\theta_{b}$ is defined in (\ref{minc00}).

%\vspace{-0.2cm}

\subsection{Feasible Regions}

We can now analyze, when other parameters (e.g., load predictions) are given, the feasible regions for $f^u,f^d,$ and $\theta$:
\begin{align}
   & 0\le f^u \le \bar{f}^u(f^d) \doteq \text{MaxUR}(\infty,f^d), \label{boundfu} \\
   & 0 \le f^d \le \bar{f}^d(f^u) \doteq \text{MaxDR}(\infty,f^u), \label{boundfd} \\
   & \text{MinC}(0,0) \le \theta \le \bar{\theta} \doteq \max \text{MinC}(f^u,f^d). \label{boundtheta}
\end{align}
The infinity argument (i.e., $\infty$) in (\ref{boundfu}) and (\ref{boundfd}) implies that the corresponding constraint has been relaxed. Thus, we can conclude that the feasible regions of $f^u$ and $f^d$ are coupled. Note that there might be multiple pairs $(f^u,f^d)$ such that DS$(f^u,f^d)=0$. We regard all such pairs as the non-interesting region \rev{for analytical purposes.} We denote the boundary of this region by $\underline{f}^u(f^d)$ and $\underline{f}^d(f^u)$, where
\begin{align}
   & \underline{f}^u(f^d) = \text{MaxUR}(0,f^d), \label{boundful} \\
   & \underline{f}^d(f^u) = \text{MaxDR}(0,f^u). \label{boundfdl}
\end{align}
Clearly, if $\underline{f}^u(f^d)=\bar{f}^u(f^d)$ and $\underline{f}^d(f^u)=\bar{f}^d(f^u)$, flexible ramping products will not impose any additional cost. Hence, for the subsequent analysis, we only concentrate on cases $(f^u,f^d,\theta)$ for which
\begin{align}
    \underline{f}^u(f^d)<& \ f^u < \bar{f}^u(f^d), \label{boundfun} \\
    \underline{f}^d(f^u)<& \ f^d < \bar{f}^d(f^u), \label{boundfd2} \\
    0<&\  \ \theta  \ < \bar{\theta}. \label{boundthetan}
\end{align}
%with $\underline{f}^u(f^d)<\bar{f}^u(f^d)$ and $\underline{f}^d(f^u)<\bar{f}^d(f^u)$.

%\vspace{-0.3cm}
\section{Analytical Understanding}
\label{analytical}

\subsection{Analytical Relationships among the Functions}
The feasible regions (\ref{boundfu})-(\ref{boundfdl}) shed light on some basic properties of the three functions on the boundaries. In this section, we provide two theorems that further ascertain the underlying relations between the three functions. First, based on the results in linear parametric programming \cite{Holder2010}, we can formulate the following theorem:

\vspace{0.2cm}
\begin{theorem}\label{thm1}
\emph{(a) The} MinC \emph{function is continuous, piecewise linear, convex, and non-decreasing in both $f^u$ and $f^d$.\\ (b) The} MaxUR \emph{function is continuous, piecewise linear, and concave in both $\theta$ and $f^d$; it is non-decreasing in $\theta$ while non-increasing in $f^d$. \\(c) The} MaxDR \emph{function is continuous, piecewise linear, and concave in both $\theta$ and $f^u$; it is non-decreasing in $\theta$ while non-increasing in $f^u$.}
\end{theorem}
\vspace{0.2cm}

Over the region of interest, defined by (\ref{boundfun})-(\ref{boundthetan}), the three functions become strictly monotone functions. Hence, given any of the two arguments, their inverse functions exist over the region of interest. We can further show the key results:

\vspace{0.2cm}
\begin{theorem}\label{thm2}
\emph{In the region (\ref{boundfun})-(\ref{boundthetan}), given any $f^u_0,f^d_0,\theta_0$, \\(a)} MinC$(f^u,f^d_0)$ \emph{and} MaxUR$(\theta,f^d_0)$ \emph{are inverse functions of each other; (b)} MinC$(f^u_0,f^d)$ \emph{and} MaxDR$(\theta,f^u_0)$ \emph{are inverse functions of each other;  (c)} MaxUR$(\theta_0,f^d)$ \emph{and} MaxDR$(\theta_0,f^u)$ \emph{are inverse functions of each other.\\ Mathematically, we have}
\begin{align}
    \text{MaxUR}(\text{MinC}(f^u,f^d),f^d) & = f^u, \\
    \text{MinC}(\text{MaxUR}(\theta,f^d),f^d) \ & = \theta, \\
    \text{MaxDR}(f^u,\text{MinC}(f^u,f^d)) & = f^d,\\
    \text{MinC}(f^u,\text{MaxDR}(\theta,f^u)) \ & = \theta, \\
    \text{MaxUR}(\theta,\text{MaxDR}(\theta,f^u)) & = f^u,\\
    \text{MaxDR}(\theta,\text{MaxUR}(\theta,f^d)) & = f^d.
\end{align}
\end{theorem}

The proofs for these  two  theorems are presented in Appendix \ref{app1} and \ref{app2}, respectively. Based on Theorem \ref{thm2}, we can concentrate our analysis on only one of the three functions, and it enjoys ``triple optimality''. For instance, for any point $(f^u_0,f^d_0,\theta_0)$ on the MinC function, it apparently means given $f^u_0$ and $f^d_0$, the minimal distortion cost is $\theta_0$. With Theorem \ref{thm2}, we can also argue that given financial budget $\theta_0$ and up-ramping requirement $f^u_0$, the maximal down flexible ramping capacity that the system can contribute is $f^d_0$; and given financial budget $\theta_0$ and down-ramping requirement $f^d_0$, the maximal up flexible ramping capacity is $f^u_0$.

\subsection{Efficient Construction}
\label{routine}

\begin{figure*}[t]
\begin{center}
\includegraphics[width=1.8\columnwidth]{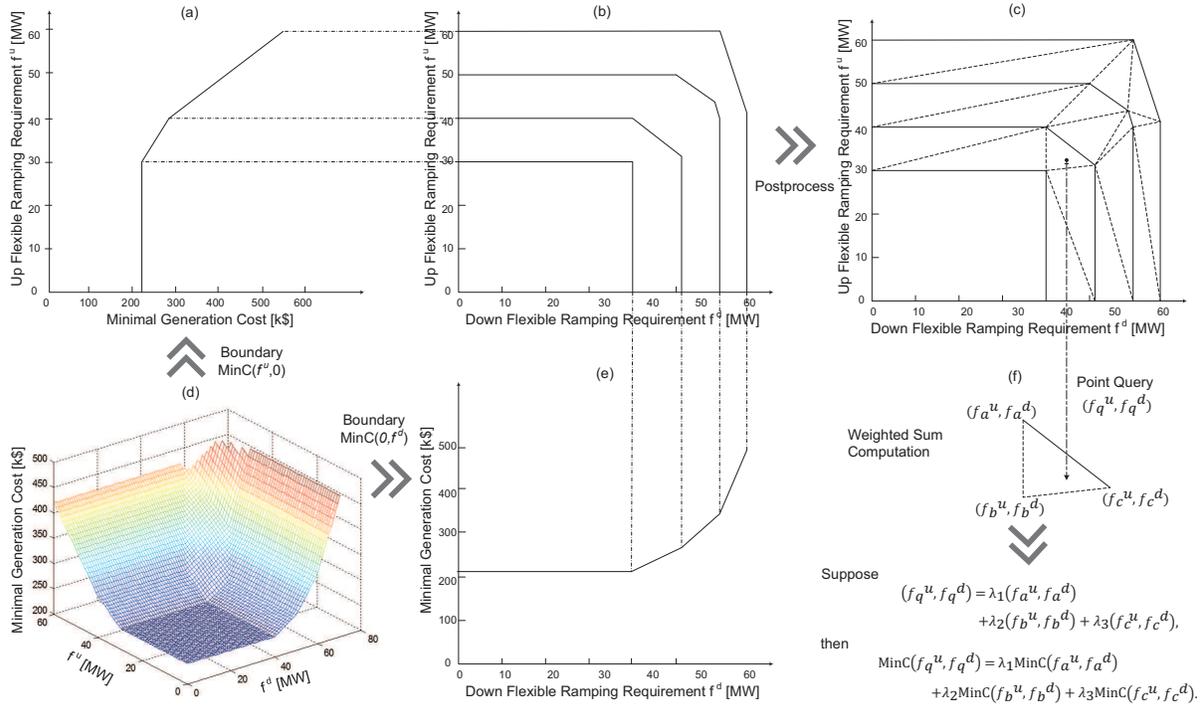}
%\vspace{-0.3cm}
 \caption{Visualization of Algorithm 2. Reading order: (d)$\rightarrow$(a),(e)$\rightarrow$(b)$\rightarrow$(c)$\rightarrow$(f). (a) Boundary of sample MinC function, MinC($f^u,0$); (b) Results returned by the Main Algorithm in Algorithm 2; (c) Results returned by the postprocess in Algorithm 2; (d) Sample MinC function to be constructed; (e) Boundary of sample MinC function, MinC($0,f^d$); (f) Example to query the value of $(f_q^u,f_q^d)$.}\vspace{-0.5cm}
 \label{fig:algo}
\end{center}
\end{figure*}

We now discuss the computational efforts required to construct the three functions. The evaluation of each function essentially corresponds to solving an OPF problem for any combination of values for the parameters of the function. Hence, even with the DC approximation, it is still computationally expensive for large power networks to determine the function over the full range of feasible parameter values \cite{purchala2005usefulness}. Hence, in practice, it may be difficult to directly compute the functions, even just for evaluation purposes.

Fortunately, the properties of the three functions can help to substantially reduce the computational efforts. Let us take the MinC function as an example: the function is piecewise linear and non-decreasing in both arguments. Thus, using Lagrangian duality \cite{boyd2004convex}, we can characterize the slopes of the piecewise linear segments and use these slopes to provide an efficient way to compute the function. If there is a single argument, \cite{wu2013unifying} gives an algorithm to construct the function with $m$ linear segments in $O(m)$ steps. \rev{Our construction of the parametric optimization functions will rely on the following single parameter function construction subroutine:}

\vspace{0.2cm}
\noindent \textbf{Algorithm 1: Single Parameter Function Construction}
\vspace{0.1cm}

\noindent Example: Construct MaxUR($\theta_i,f^d$) for given $\theta_i$ in $[a,b]$.
\begin{enumerate} \itemsep 2pt
  \item Compute MaxUR($\theta_i,a$) and MaxUR($\theta_i,b$), and obtain the Lagrangian multipliers associated with constraint $\sum_{n\in\mathcal{N}} r_{n,1}^d = f^d$ for the two cases (denoted by $\lambda_a$ and $\lambda_b$, respectively). Solve the following system of equations:
      \begin{align}\label{systemofeqs}
        & c_{a,b} - \text{MaxUR}(\theta_i,a) = \lambda_a (f_{a,b}^d - a),\\
        & c_{a,b} - \text{MaxUR}(\theta_i,b) = \lambda_b (f_{a,b}^d - b),
      \end{align}
      to obtain ($f_{a,b}^d,c_{a,b}$), where $a\le f_{a,b}^d\le b$.
  \item If MaxUR$(\theta_i,f_{a,b}^d)=c_{a,b}$, then within interval [$a,b$]:
  \begin{equation}\label{MaxURconstruction}
  \begin{aligned}
   \!\!\!\!\!\! \!\!\!\!\!\!& & & \text{MaxUR}(\theta_i,f^d) \\
   \!\!\!\!\!\! \!\!\!\!\!\!& & =& \begin{cases}
    \text{MaxUR}(\theta_i,a)\!+\!\lambda_a (f_{a,b}^d \!\!-\! a),  \text{if }a\!\le\! f^d\!\le\! f_{a,b}^d,\!\\
    \text{MaxUR}(\theta_i,b)\!+\!\lambda_b (f_{a,b}^d \!-\! b),  \text{if }f_{a,b}^d\!\le\! f^d\!\le\! b.
    \end{cases}
    \end{aligned}
  \end{equation}
  If MaxUR$(\theta_i,f_{a,b}^d) \ne c_{a,b}$, then construct the function over intervals [$a,f_{a,b}^d$] and [$f_{a,b}^d,b$], respectively.
\end{enumerate}

\vspace{0.2cm}
\rev{This algorithm will return all the breaking points of the function as well as  the slope of each segment. The breaking point is the intersection point of two adjacent line segments in the piecewise linear function. This algorithm can be easily applied to construct MinC($f_0^u,f^d$) for given $f^u_0$, and MinC($f^u,f_0^d$) for given $f^d_0$. The corresponding Lagrangian multipliers for constructing MinC($f_0^u,f^d$) are those associated with constraint $\sum_{n\in\mathcal{N}} r_{n,1}^d = f^d$, while the Lagrangian multipliers for constructing MinC($f^u,f_0^d$) are those associated with constraint $\sum_{n\in\mathcal{N}} r_{n,1}^u = f^u$. In practice, the Lagrangian multipliers can be obtained by a variety of primal-dual solvers for convex optimization problems, such as CVX \cite{cvx,gb08}.
\vspace{0.1cm}

For the single argument function construction, it suffices to utilize the piecewise linearity and monotonicity as shown in Algorithm 1. However, for the two-argument functions discussed in this paper, these two properties are not enough. By utilizing the additional triple optimality property, we propose two routines to construct the functions: one for computation, and the other for visualization. Both routines will use Algorithm 1 as a subroutine.

Note that the MinC function is piecewise linear in both arguments. Therefore, it comprises several facets. The general idea to construct the MinC function is to efficiently identify the breaking lines between different facets (i.e., the boundaries of the facets). The horizontal section of a fixed $\theta$ is MaxUR($\theta,f^d$) function. As $\theta$ increases, the number of segments in the MaxUR($\theta,f^d$) function may change. Any such change corresponds to an emerging facet, or a vanishing one. Such changes will also be reflected on the boundaries, specifically, MinC($0,f^d$) and MinC($f^u,0$), as breaking points on one or both boundaries. Note that, the breaking points are given by Algorithm 1. Therefore, it suffices to track all the breaking points on both boundaries to identify the breaking lines between different facets.

\vspace{0.2cm}
\noindent \textbf{Algorithm 2: MinC Construction for Computation}
\vspace{0.1cm}

\noindent \textbf{Preprocess:}

Construct the boundaries: MinC($0,f^d$), and MinC($f^u,0$) (Algorithm 1).
Denote the breaking points in MinC($0,f^d$) by $(f_1^d,\cdots,f_s^d)$, and the corresponding costs by $(\theta_1^d,\cdots,\theta_s^d)$. Denote the breaking points in MinC($f^u,0$) by  $(f_1^u,\cdots,f_v^u)$, and the corresponding costs by $(\theta_1^u,\cdots,\theta_v^u)$.
%\noindent Example: Construct MaxUR($\theta_i,f^d$) for given $\theta_i$ in $[a,b]$.
%  \item Rank $(\theta_1^d,\cdots,\theta_s^d)$ and $(\theta_1^u,\cdots,\theta_v^u)$ in an ascending order. Denote the rank by $(\theta_1,\cdots,\theta_{s+v})$.
%  \item For $i$ from $1$ to $s+v$, if $\theta_i=\theta_{i+1}$, denote
%  \item Suppose $s\le v$. Then, for $i$ from $1$ to $s$, construct MaxUR($\theta_i,f^d$), where
%  \begin{equation}\label{cons1}
%    \theta_i = \text{MinC}(f_i^u).
%  \end{equation}
%  Note that the two end points of  MaxUR($\theta_i,f^d$) are $(0,f^d_i,\theta_i)$ and $(f^u_i,0,\theta_i)$.
%  \item If $s=v$, then skip Step 5.
%  \item For $i$ from $s+1$ to $v$, we continue constructing the MaxUR functions. The only difference is that now we only know one end point $(f_i^u,0,\theta_i)$, where $\theta_i$ still can be computed by (\ref{cons1}). The other end point can be given by $(f_s^d,f_i^{u,\star},\theta_i)$, where
%      \begin{equation}\label{cons2}
%        f_i^{u,\star} = \text{MaxUR}(\theta_i,f_s^d).
%      \end{equation}
%      Then, construct MaxUR($\theta_i,f^d$) over $[f_i^{u,\star},f_i^d]$ for the fixed $\theta_i$.
%  \item

\vspace{0.1cm}
\noindent \textbf{Main Algorithm:}

\begin{algorithmic}
 \STATE $i \leftarrow 1$; $j\leftarrow 1$;
 \WHILE{$i\le s$ and $j\le v$}
 \IF{$\theta_i^d = \theta_j^u$}
 \STATE Construct MaxUR($\theta_i^d,f^d$) over $[0,f^d_i]$;
 \STATE $i \leftarrow i+1$; $j\leftarrow j+1$;
 \ELSIF{$\theta_i^d < \theta_j^u$}
 \STATE Construct MaxUR($\theta_i^d,f^d$) over $[0,f^d_i]$;
 \STATE $i \leftarrow i+1$;
 \ELSE
 \STATE Construct MaxDR($\theta_j^u,f^u$) over $[0,f^u_j]$;
 \STATE $j \leftarrow j+1$;
 \ENDIF
 \ENDWHILE
 \end{algorithmic}

\vspace{0.1cm}
\noindent \textbf{Postprocess:}

Now we obtain the horizonal sections where there is either an emerging facet or a vanishing one. According to the convex property of the MinC function, we can connect the adjacent horizonal sections and partition the feasible region (\ref{boundfun})-(\ref{boundthetan}) into triangles, such that there is no breaking line in any of the triangles. Fig. \ref{fig:algo} visualizes this algorithm.

\vspace{0.2cm}

The final step returns a set of triangles. For any given parameters $(f^u,f^d)$, one can query which triangle it belongs to with a binary search (or more advanced query techniques, see \cite{vornoi} for more details). Then, the value of MinC$(f^u,f^d)$ is given by the weighted sum of the three end points of the triangle. An efficient query is crucial for our subsequent linear search algorithm to solve the risk-limiting economic dispatch problem.

We want to emphasize that the number of triangles is limited. Suppose we have $\alpha$ decision variables and $\beta$ inequality constraints. Even in the worst case, instead of having $2^{\beta+\alpha}$ triangles, there will be at most $O(\beta \alpha^{1/3})$ triangles \cite{dey1997improved}. Theoretically, this is already very efficient since the number is almost linear in the input size. In practice, the MinC function is likely to be partitioned into only dozens of triangles (as shown in the case studies), and hence is very efficient to construct.}

%\vspace{0.1cm}

\begin{figure*}[t]
\begin{center}
\includegraphics[width=1.7\columnwidth]{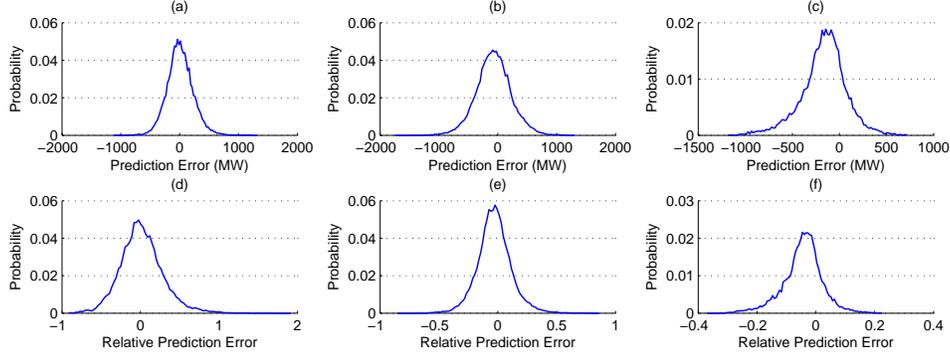}
%\vspace{-0.3cm}
 \caption{Renewable energy prediction error model. Prediction error distribution for wind power output (a) between 10\%-30\%; (b) between 30\%-70\%; (c) above 70\%. Relative prediction error distribution for wind power output (d) between 10\%-30\%; (e) between 30\%-70\%; (f) above 70\%.}\vspace{-0.3cm}
 \label{fig:errordis}
\end{center}
\end{figure*}

Although Algorithm 2 is sufficient for computation, it does not provide an intuitive way to visualize the functions. Therefore, we devise an efficient way to obtain the contour of the MinC function with $k$ lines in Algorithm 3:

\vspace{0.2cm}
\noindent \textbf{Algorithm 3: MinC Construction for Visualization}
\vspace{0.1cm}

\begin{enumerate} \itemsep 2pt
  \item Calculate $\bar{\theta}$ by solving the following problem:
  \begin{equation}\label{defthetabound}
    \begin{aligned}
  & &  \bar{\theta}= \underset{f^u,f^d}{\max} & \ \ \text{MinC}(f^u,f^d) \\
   & &  s.t. & \ \ 0\le f^u \le \text{MaxUR}(\infty,0),\\
   & & & \ \ 0\le f^d \le \text{MaxDR}(\infty,0).
    \end{aligned}
 \end{equation}
  \item To obtain the contour of the MinC function, divide the interval [MinC$(0,0), \bar{\theta}$] into $k-1$ equally incremental sub-intervals to draw the contour with $k$ lines. Denote the $k$ end points by $\theta_1,\cdots,\theta_k$.
  \item The contour line with the same cost $\theta_i$ is simply a MaxUR (or equivalently, MaxDR) function with fixed cost $\theta_i$. We may again refer to Algorithm 1 to construct MaxUR($\theta_i,f^d$) in the interval  [$0,\text{MaxDR}(\theta_i,0)$].
\end{enumerate}

\vspace{0.2cm}
Using this routine to construct a MinC function with $k$ lines, $O(\sum_{i=1}^k m_i)$ optimization problems need to be solved, where $m_i$ is the number of line segments of function MaxUR($\theta_i,f^d$) with fixed $\theta_i$. \rev{Based on the similar argument in the Algorithm 2 analysis, all the $m_i$'s are also almost linear in the input size of the problem.}

\section{Risk-limiting Energy Dispatch}
\label{distortion}

Bearing the  relationship between the generation cost and the ramping capacity requirements in mind, we now seek to understand the other dependency - how the ramping capacity requirements depend on the risk-limiting constraint, i.e., to solve the risk-limiting economic dispatch problem. Towards exploiting this dependency, in this section, we first discuss the renewable energy prediction error model. Then, we introduce the linear search algorithm to obtain the optimal combination of the key parameters for flexible ramping products.

\subsection{Renewable Energy Prediction Error Model}
\label{errorModelSection}

We use the Bonneville Power Administration (BPA) predicted and actual wind power data \cite{BPAdata} \rev{with a temporal resolution of 5 minutes} to obtain the prediction error model. For a historic dataset for a wind plant with the maximal capacity of 4,500 MW in BPA, we first note that when the power output is less than 10\% of the maximal capacity, the relative prediction error can be extremely large (or even arbitrarily large when the actual power output is zero) while the amount of the prediction error is relatively small, which implies the flexible ramping requirements in this case are not critical. Therefore, we trim these data from the dataset. Then, we divide the trimmed historical data into three groups: 10\%-30\% of the maximal capacity, 30\%-70\% of the maximal capacity, and above 70\% of the maximal capacity. The (relative) prediction error distributions for these three groups are illustrated in Fig. \ref{fig:errordis}. Since prediction errors need to be compensated by ramping generation, the level of prediction errors determines how much ramping capacity is required.

Consequently, we can compare the three cases in Fig. \ref{fig:errordis}. When the wind power output is between 10\%-30\% of the capacity, the mean prediction error is 7.8 MW with a standard deviation of 223 MW. The relative prediction error is also quite significant in this case, with a mean value of 0.015 and a standard deviation of 0.27. With the increase of the wind power output (between 30\%-70\% of the capacity, which is the most common range of power output of a wind plant), the relative prediction error drops significantly (with a standard deviation of only 0.14), but the mean value of prediction error shifts to -70.7 MW and its standard deviation is 300 MW. On windy days, with wind power output of more than 70\% of the capacity, since we know the maximal capacity of the wind plant, the prediction seems to perform reasonably well, with a mean value of -77.7 MW and a standard deviation of only 175 MW. Also, thanks to this upper bound, the standard deviation of the relative prediction error is now only 0.05.

%The relative prediction error decreases with the increase of the wind power generation, though the range of absolute prediction error increases.

\subsection{Risk-limiting Energy Dispatch}
\label{minidispatch}

Assuming prediction error distributions as shown in Fig. \ref{fig:errordis}, we seek to achieve the minimal generation cost that satisfies the risk-limiting constraints. Mathematically, if a confidence level of $p$\% is desired, the ISO needs to solve the following optimization problem:
\begin{equation}\label{min_DS}
\begin{aligned}
  & &   \underset{f^u,f^d}{\min} & \ \ \text{MinC}(f^u,f^d) \\
   & &  s.t. & \ \ \Pr ( \boldsymbol{1}^T \boldsymbol{\hat{d}}_1 \!-\! f^d \! \le \!\boldsymbol{1}^T \boldsymbol{d}_1\! \le\! \boldsymbol{1}^T \boldsymbol{\hat{d}}_1\!+\!f^u ) \ge p\%
    \end{aligned}
\end{equation}

Since we use the actual prediction error distributions, there are no symmetric nor other nice analytical properties. Therefore, we propose  a linear search method to obtain the minimal generation cost. \rev{In particular, for any given parameters ($f^u,f^d$), instead of solving the OPF problem (\ref{classical_energy_dispatch}), MinC($f^u,f^d$) can be efficiently obtained by querying the set of triangles returned by Algorithm 2.

To formally describe the linear search algorithm, we first divide the search region $[a,b]$ (given by the prediction error probability distribution) into intervals, each being of length $\delta$. Suppose $a<0$ and $b>0$. Then, we employ two for-loops to enumerate all the possible combinations $(f^d_0,f^u_0)$ that satisfies the risk-limiting constraint.

For any given $f^d_0$ specified by the outer for-loop, the inner for-loop tries to identify the shortest interval $[-f^d_0,f^u]$, which satisfies the risk-limiting constraint. Suppose $f^u_0$ is the desired parameter to form the shortest interval. Then, we update $f^u_{start}$ with $f^u_0$, and conduct the comparison to see if the new combination $(f^d_0,f^u_0)$ achieves a lower generation cost. Due to the monotonicity of MinC function, there is no need to query MinC($f^u,f^d_0$), for any $f^u>f^u_0$. Thus, we break the inner for-loop. After that, the outer loop starts again, and sets $f^d$ to be $f^d+\delta$. Based on the monotonicity of the cumulative probability distribution, the inner for-loop search process can directly start from the updated $f^u_{start}$.

\vspace{0.2cm}
\noindent \textbf{Algorithm 4: Linear Search Optimal Combination}
\vspace{0.1cm}

 \begin{algorithmic}
 \STATE $f^u_{start} \leftarrow 0$; $f^u_{opt}\leftarrow 0$; $f^d_{opt}\leftarrow 0$; $opt \leftarrow \infty$;
 \FOR{$f^d = -a:\delta:0$}
 \FOR{$f^u = f^u_{start}:\delta:b$}
 \IF{$(f^u,f^d)$ satisfies the risk-limiting constraint}
  \IF{MinC$(f^u,f^d) \le opt$ }
 \STATE $opt \leftarrow \text{MinC}(f^u,f^d)$
 \STATE $f^u_{opt} \leftarrow f^u$
 \STATE $f^d_{opt} \leftarrow f^d$
 \ENDIF
 \STATE $f^u_{start} \leftarrow f^u$
 \STATE Break;
 \ENDIF
 \ENDFOR
 \ENDFOR
 \end{algorithmic}
 \vspace{0.2cm}

This search process will query MinC($f^u,f^d$) at most $(b-a)/\delta$ times, which is linear in the length of the interval. However, we want to emphasize that although the linear search is efficient, its accuracy relies on the selection of $\delta$. If the problem were convex, we could have implemented binary search to achieve arbitrary accuracy. From another point of view, this is also the evidence for the hardness of the non-convexity.}

It is worth noting that, the error distributions might be approximated by Gaussian or other well-studied distributions, which can lead to improved accuracy in the solution of the optimization problem but will lead to errors introduced by the approximation of the probability distributions.
%----------------------------------------------------------------------
% SECTION VI: Simulation Results
%----------------------------------------------------------------------

\section{Case Study}
\label{simulation}

\begin{figure*}[t]
\begin{center}
\subfigure[3-bus system.] {
\includegraphics[width=0.32\columnwidth]{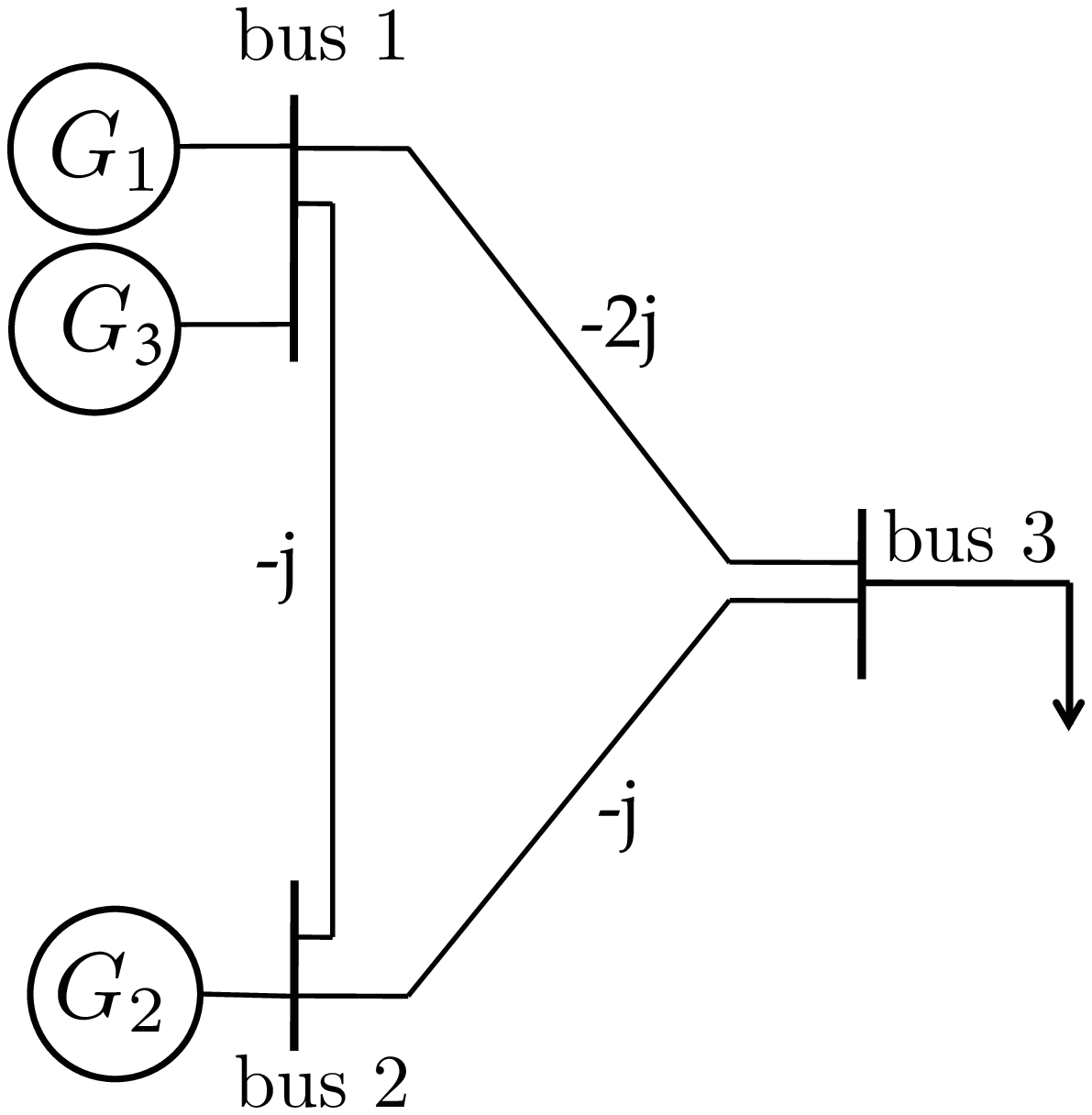}
}
\subfigure[DS with fixed $f^d$.] {
\includegraphics[width=0.36\columnwidth]{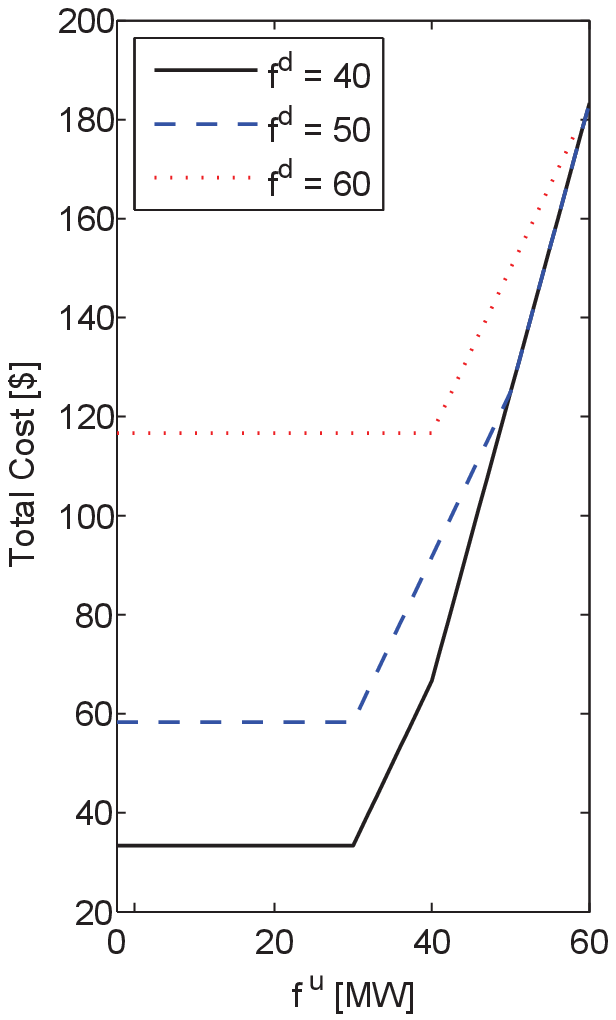}
}
\subfigure[DS with fixed $f^u$.] {
\includegraphics[width=0.35\columnwidth]{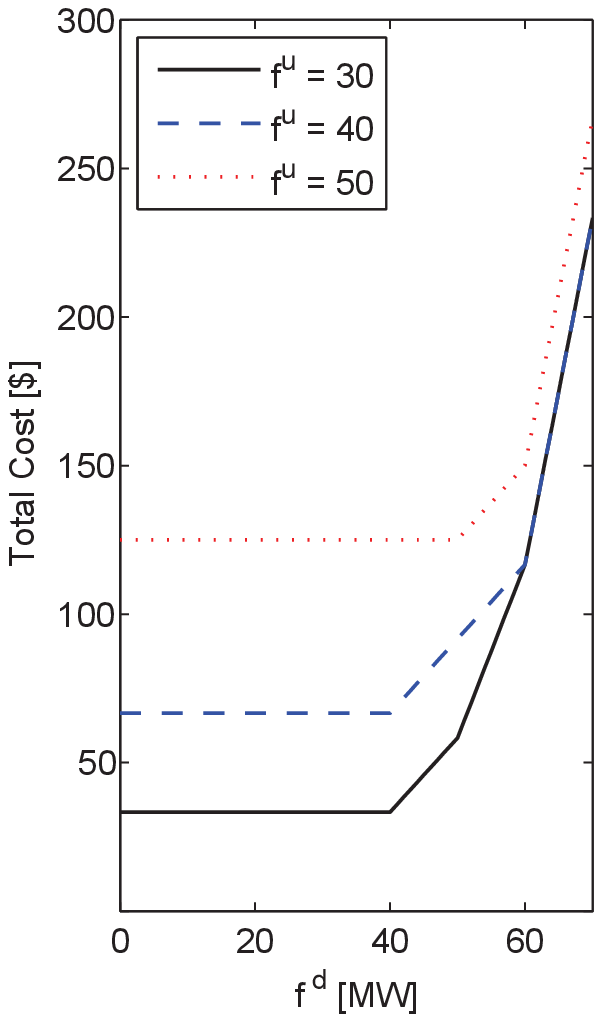}
}
\subfigure[Contour of the DS function.] {
\includegraphics[width=0.75\columnwidth]{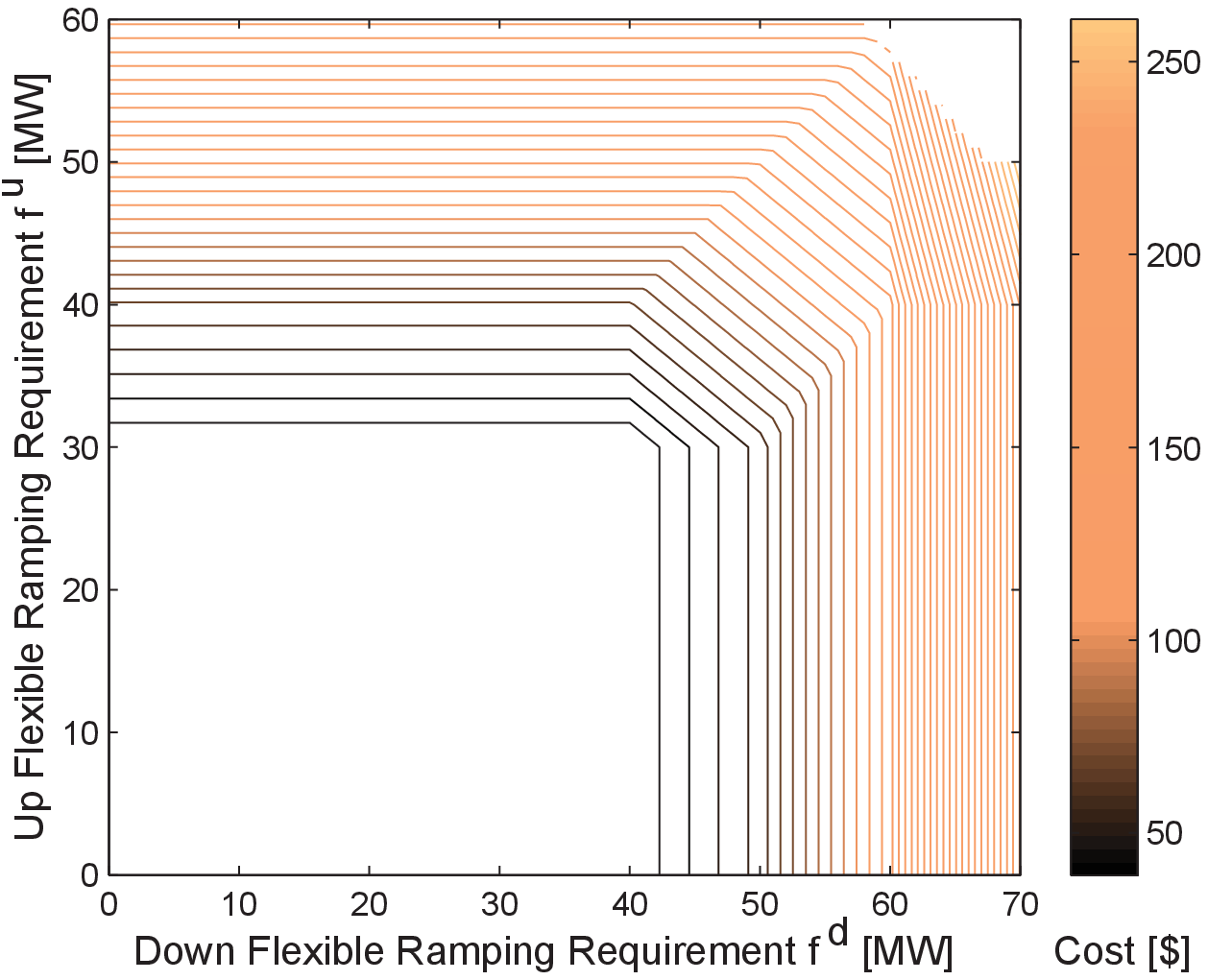}
}
 \caption{System information and simulation results for the 3-bus system.}%\vspace{-0.3cm}
 \label{fig:3-bus-wo-line}
\end{center}
\end{figure*}

In this section, we first consider a prototype 3-bus system and Garver's 6-bus system to highlight the properties of our parametric functional approach. Then, we turn to  more realistic scenarios by evaluating the influence of flexible ramping products on the WECC 240-bus system. Both simulation results reveal interesting information on the relationship between the generation cost and the key parameters of flexible ramping products. We hope such information can help ISOs to better evaluate the new products and develop methods to achieve the most cost effective market for the new products.

\rev{To better quantify the distortion incurred by the new products, we define the distortion cost function DS($f^u,f^d$) as the difference in generation cost between the optimal economic dispatch with flexible ramping requirements ($f^u,f^d$) and the optimal dispatch without flexible ramping requirements, i.e.,
\begin{equation}\label{DS_def}
    \text{DS}(f^u,f^d) = \text{MinC}(f^u,f^d)-\text{MinC}(0,0).
\end{equation}
%Since the quantity MinC(0,0) will play an important role throughout the paper, for notational simplicity, we denote
%\begin{equation}\label{minc00}
%    \theta_b = \text{MinC}(0,0).
%\end{equation}
%Note that MinC(0,0) corresponds to the conventional economic dispatch problem.

Note that DS function is simply the shifted MinC function. Therefore, it preserves all the analytical properties of MinC function. For the subsequent analysis, we will demonstrate the performance of our approach with DS functions.}

\subsection{Prototype 3-bus System}

\begin{table}[b]
\caption{Generator Information for the 3-bus System.}\vspace{-0.3cm} \label{Table:info}
\begin{center}
\begin{tabular}{|c|c|c|c|c|}
\hline
& $C_n$ [\$/MW] & $\Delta g_n$ [MW/5 min] & ${g}_{n,-\!1}$ [MW] &$\bar{g}_n$ [MW]  \\
\hline
$G_1$ & 50 & 20 & 90 & 100 \\
$G_2$ & 120 & 30 & 0& 100 \\
$G_3$ & 80 & 20 & 20& 20  \\
\hline
\end{tabular}
\end{center}\vspace{-0.3cm}
\end{table}

We illustrate Algorithm 2 and Algorithm 3 using a prototype 3-bus system. As shown in Fig. \ref{fig:3-bus-wo-line}(a), there are two generators $G_1$ and $G_3$ at bus 1; the third generator $G_2$ is at bus 2; and bus 3 is a load bus. Table \ref{Table:info} gives the necessary data for the three generators. We assume $\underline{g}_n$'s are all zero for this prototype 3-bus system. We do not consider the line capacity constraints in this example. Assume the generators are to be dispatched for a forecasted net load of 110 MW at $t=0$, and 120 MW at $t=1$. Then, without any flexible ramping requirement, the optimal economic dispatch profile at $t=0$ is (100, 0, 10) MW, while it is (100, 0, 20) MW at $t=1$.

Note that, such economic dispatch comes with free ramping capacity. To reserve ramping capacities for time $t=1$, the free up-ramping capacity is 30 MW from $G_2$, which corresponds to the horizontal segment in the lower envelope of Fig. \ref{fig:3-bus-wo-line}(b). We want to highlight that such horizontal region is precisely the non-interesting region discussed in Section \ref{functional}. After this free and non-interesting region, the two parameters start distorting the generation output profile compared with the optimum given by DS$(0,0)$. Similarly, we can analyze the free down-ramping capacity for $t=1$: 20 MW from $G_1$ plus 20 MW from $G_3$, which corresponds to the horizontal segment in the lower envelope of Fig. \ref{fig:3-bus-wo-line}(c). Fig. \ref{fig:3-bus-wo-line}(b) and (c), together with the contour of the DS function shown in Fig. \ref{fig:3-bus-wo-line}(d), illustrate all the properties (piece-wise linearity, monotonicity, and convexity) stated in Theorem \ref{thm1}.

Following Algorithm 3 in Section \ref{routine}, to efficiently construct the DS function for visualization, we first identify $\bar{\theta}$, which according to Fig. \ref{fig:3-bus-wo-line}(d) is given by DS(50,70). We divide [$0,\bar{\theta}$] into 29 slots to obtain the contour of the DS function with 30 lines. In this example, from the constructed function shown in Fig. \ref{fig:3-bus-wo-line}(d), each contour line consists of at most three segments. Hence, the routine requires solving at most 5 optimization problems to construct each contour line.

\begin{figure}[t]
\begin{center}
 \includegraphics[width=7cm]{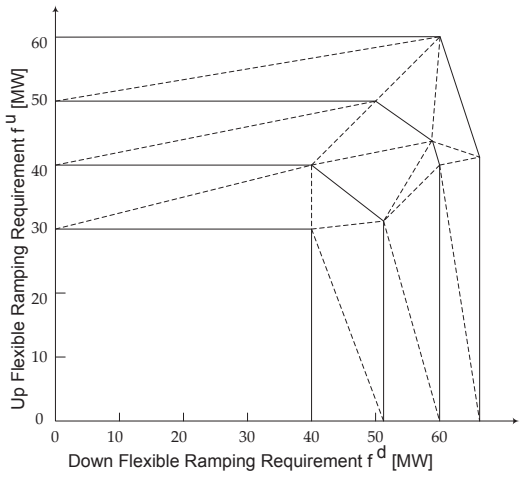}
 %\vspace{-0.3cm}
 \caption{Triangulation based on DS function.}%\vspace{-0.3cm}
 \label{fig:3bus_triangle}
\end{center}
\end{figure}

\rev{To efficiently construct the function for computational purposes, we can follow Algorithm 2 to obtain the set of triangles as shown in Fig. \ref{fig:3bus_triangle}. The solid lines are given by the main algorithm in Algorithm 2. Based on these solid lines, and the convexity of MinC function, we can conduct the postprocessing part of the algorithm to obtain all the triangles (the dashed lines). Note that there could be different sets of triangles to partition the space. However, the total numbers of triangles for different sets are the same. The system model has 9 decision variables and 27 constraints. Instead of having $O(2^{36})$ triangles, there are altogether 29 triangles. This also confirms that the number of triangles will be almost linear in the input size. In fact, as the system scales up, the number of triangles does not grow linearly in practice in that there are often limited binding constraints.}

\subsection{Garver's 6-bus System}

\begin{figure}[t]
\begin{center}
 \includegraphics[width=4.5cm]{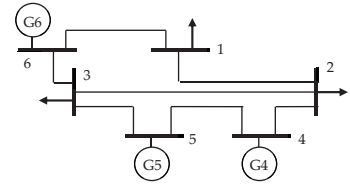}
 %\vspace{-0.3cm}
 \caption{Garver's 6-bus system.}%\vspace{-0.3cm}
 \label{fig:6bus}
\end{center}
\end{figure}

\begin{figure}[t]
\begin{center}
 \includegraphics[width=7.5cm]{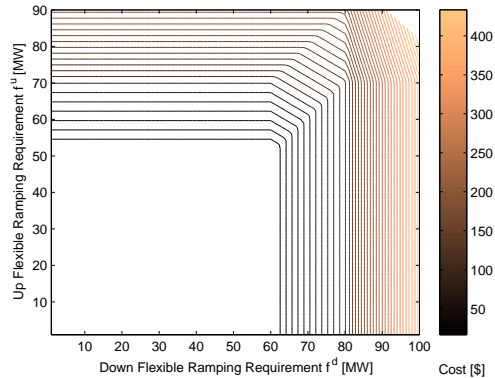}
% \vspace{-0.3cm}
 \caption{Contour of the DS function for Garver's 6-bus system.}%\vspace{-0.1cm}
 \label{fig:contour_6bus}
\end{center}
\end{figure}

\begin{table}[b]
\caption{Generator Information for Garver's 6-bus System.}\vspace{-0.3cm} \label{Table:info_6bus}
\begin{center}
\begin{tabular}{|c|c|c|c|c|}
\hline
& $C_n$ [\$/MW] & $\Delta g_n$ [MW/5 min] & ${g}_{n,-\!1}$ [MW] & $\bar{g}_n$ [MW] \\
\hline
$G_4$ & 58 & 50 & 0& 200  \\
$G_5$ & 52 & 20 & 130&  150 \\
$G_6$ & 54 & 20 & 30 &40 \\
\hline
\end{tabular}
\end{center}\vspace{-0.2cm}
\end{table}

We analyze Garver's 6-bus system \cite{4074249} shown in Fig. \ref{fig:6bus}. Table \ref{Table:info_6bus} provides the necessary information for the generators at buses 4, 5, and 6. Again, we assume $\underline{g}_n$'s are all zero for this system. The initial generation outputs are 0, 130, and 30 MW, respectively. The loads are located at buses 1, 2, and 3. At time $t=0$, the forecasted net loads are 59.5 MW each, while at time $t=1$, the forecasted net loads are 63 MW each. Fig. \ref{fig:contour_6bus} shows the contour of the DS function for this case. %The maximal and minimal possible distortion costs, shown in Fig. \ref{fig:penetr_6_bus}, are smaller compared to the 3-bus case. This is largely because in this case, the generator's marginal costs are roughly the same. This raises an interesting observation: if the marginal cost of a natural gas turbine can be further reduced, then hopefully, the flexible ramping products will not introduce too much distortion cost.

%Based on the DS function, we would like to perform the minimal distortion cost economic dispatch at a certain confidence level.
Suppose now we want to achieve a 30\% renewable energy penetration level, as planned by the CAISO for the year of 2020 \cite{cal33}. And in this test system, the conventional generators contribute around 200 MW, which implies that the wind plant needs to supply 100 MW on average. Based on the current wind power technology \cite{boccard2009capacity}, the typical wind plant supplies only 20\%-40\% of its maximal capacity on average. Therefore, we assume there is a wind plant with a maximal capacity of 500 MW in the system.

\begin{figure}[t]
\begin{center}
 \includegraphics[width=8.5cm]{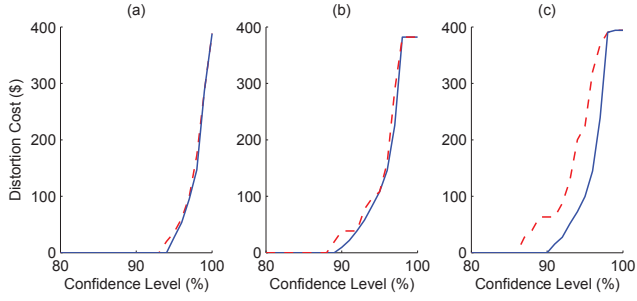}
 %\vspace{-0.3cm}
 \caption{The possible distortion costs incurred by flexible ramping products for Garver's 6-bus system: (a) 110 MW wind power; (b) 250 MW wind power; (c) 450 MW wind power. (Solid line: minimal distortion cost; dashed line: cost of the greedy approach.)}%\vspace{-0.2cm}
 \label{fig:cost_6_bus}
\end{center}
\end{figure}

Fig. \ref{fig:cost_6_bus} shows the possible distortion costs induced by flexible ramping products for the three cases discussed in Section \ref{errorModelSection} (low wind, modest wind, and high wind). \rev{We compare our risk-limiting economic dispatch approach with a greedy one. The greedy approach selects the parameters ($f^u,f^d$), which is the shortest interval (or minimal $f^u+f^d$) to achieve the risk-limiting constraint. Mathematically, ($f^u,f^d$) employed by the greedy approach is the solution to the following optimization problem:
\begin{equation}\label{greedy}
\begin{aligned}
  & &   \underset{f^u,f^d}{\min} & \ \ f^u+f^d \\
   & &  s.t. & \ \ \Pr ( \boldsymbol{1}^T \boldsymbol{\hat{d}}_1 \!-\! f^d \! \le \!\boldsymbol{1}^T \boldsymbol{d}_1\! \le\! \boldsymbol{1}^T \boldsymbol{\hat{d}}_1\!+\!f^u ) \ge p\%
    \end{aligned}
\end{equation}
%\rev{This greedy approach utilizes the near symmetric properties of the prediction error model and the DS functions. If both of them are symmetric, then this greedy approach achieves the same performance as our approach. But it is not always the case (e.g., see the mean values of the prediction error models in Section \ref{errorModelSection}).}

As shown in Fig. \ref{fig:cost_6_bus}, in all the three cases, both our approach and the greedy approach work reasonably well. Yet, our approach outperforms the greedy approach with respect to distortion cost, particularly for the cases with higher wind power output levels (Fig. \ref{fig:cost_6_bus} (b) and (c)). Numerically, compared to the greedy approach, when the confidence level is greater than 90\%, the average savings are 15.6\% when the predicted wind power is 110 MW, 21.3\% when the predicted wind power is 250 MW, and 51.3\% when the predicted wind power is 450 MW. We want to emphasize that such savings do not come with expensive computational efforts due to the piecewise linearity of the DS function as we discussed in Section \ref{minidispatch}.}

\subsection{WECC System}

\begin{figure}[t]
\begin{center}
 \includegraphics[width=7.5cm]{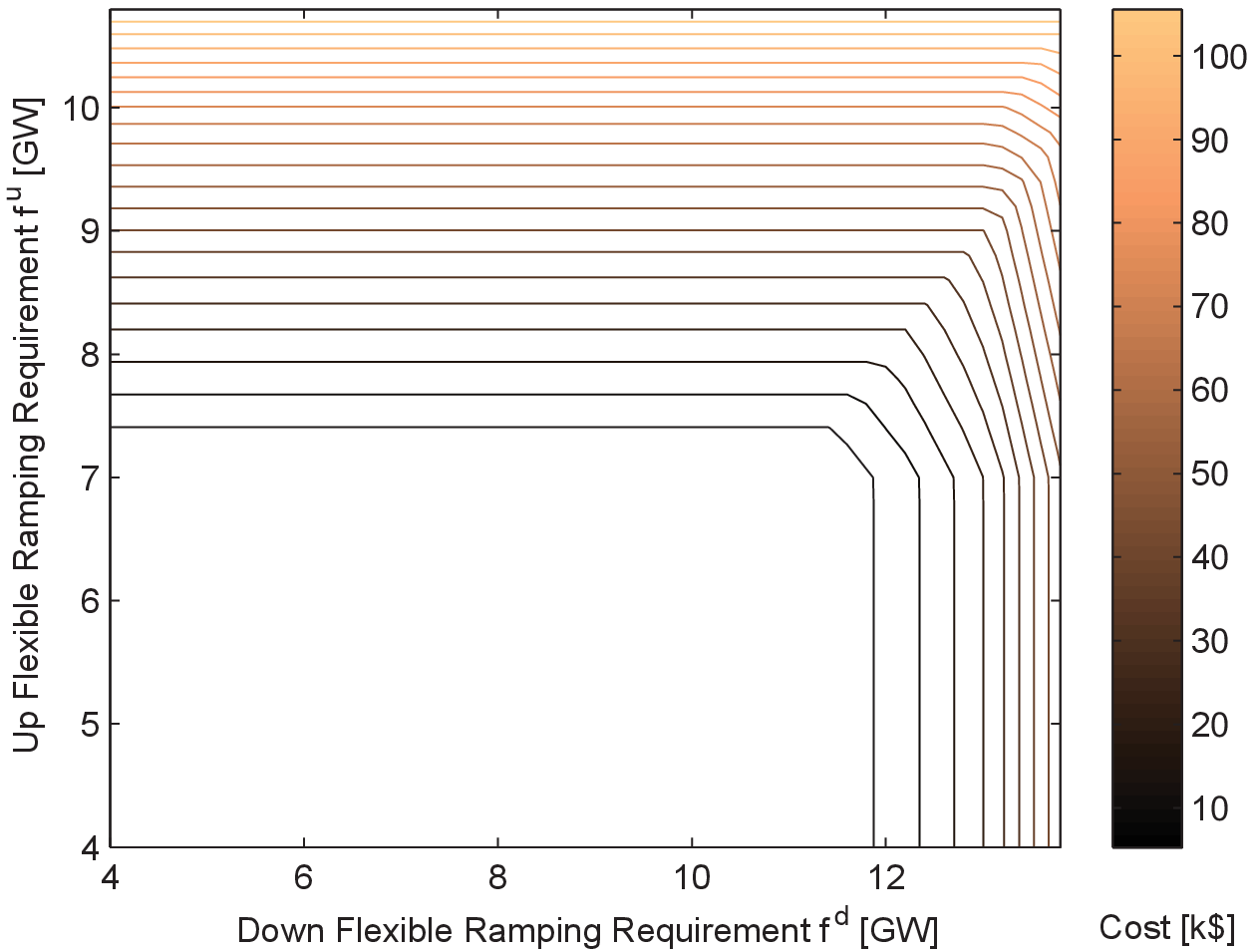}
 %\vspace{-0.3cm}
 \caption{Contour of the DS function for WECC system.}%\vspace{-0.2cm}
 \label{fig:WECC_DS}
\end{center}
\end{figure}

\begin{figure}[t]
\begin{center}
 \includegraphics[width=8.5cm]{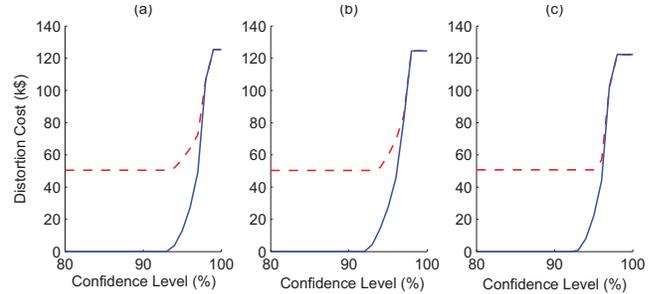}
 %\vspace{-0.3cm}
 \caption{The possible distortion costs incurred by flexible ramping products for the WECC system: (a) 16 GW wind power; (b) 35 GW wind power; (c) 70 GW wind power. (Solid line: minimal distortion cost; dashed line: cost of the greedy approach.)}%\vspace{-0.2cm}
 \label{fig:WECC_all}
\end{center}
\end{figure}

%\begin{figure}[t]
%\begin{center}
% \includegraphics[width=6.5cm]{39bus}
% \caption{IEEE 39-bus system.}\vspace{-0.3cm}
% \label{fig:39-bus-system}
%\end{center}
%\end{figure}
%
%\begin{figure}[t]
%\begin{center}
% \includegraphics[width=6.5cm]{39-dc}
% \caption{Contour of DC-based DS function for the 39-bus system.}\vspace{-0.3cm}
% \label{fig:39-dc}
%\end{center}
%\end{figure}
%
%\begin{figure}[t]
%\begin{center}
% \includegraphics[width=6.5cm]{39-ac}
% \caption{Contour of AC-based DS function for the 39-bus system.}\vspace{-0.3cm}
% \label{fig:39-ac}
%\end{center}
%\end{figure}

\begin{table}[t]
\caption{Generator Information for the WECC System.}\vspace{-0.3cm} \label{Table:info_39bus}
\begin{center}
\begin{tabular}{|c|c|c|c|c|}
\hline
Bus & $C_n$ [\$/MW] & $\Delta g_n$ [MW/5 min] &${g}_{n,-\!1}$ [MW] & $\bar{g}_n$ [MW] \\
\hline
1034 & 30 & 349.5 & 1200 & 1400 \\
1232 & 20 & 315 & 1256 & 1256 \\
1331 & 45 & 610.5 & 510 & 2438 \\
2130 & 20 & 126 & 500 & 500 \\
2637 & 55 & 27 & 30 & 110\\
4031 & 50 & 244.5 & 220 & 978 \\
4035 & 20 & 918 & 3671 & 3671 \\
4039 & 40 & 730.5 & 1000 & 2918 \\
4131 & 20 & 3240 & 12963& 12963 \\
4132 & 30 & 1422 & 4000 & 5693 \\
4231 & 20 & 904.5 & 3615 & 3615 \\
4232 & 60 & 142.5 & 150 & 573 \\
5031 & 30 & 2200.5 & 4000 & 8798.8 \\
5032 & 20 & 1102.5 & 4410.2 & 4410.5 \\
6132 & 65 & 535.5 & 550 & 2144 \\
6235 & 40 & 172.5 & 300 & 693 \\
6335 & 40 & 202.5 & 600 & 812 \\
6533 & 80 & 10.5 & 12 & 39 \\
7032 & 30 & 172.5 & 500& 691 \\
8033 & 20 & 393 & 1572 & 1572 \\
8034 & 70 & 172.5 & 180 & 688 \\
\hline
\end{tabular}
\end{center}%\vspace{-0.5cm}
\end{table}

%\begin{figure}[t]
%\begin{center}
% \includegraphics[width=7cm]{39-cost}
% \vspace{-0.2cm}
% \caption{The possible distortion costs incurred by flexible ramping products for 39-bus system: (a) 2,000 MW wind power; (b) 3,200 MW wind power; (c) 5,000 MW wind power. (Solid line: maximal distortion cost; dashed line: minimal distortion cost.)}\vspace{-0.5cm}
% \label{fig:39-cost}
%\end{center}
%\end{figure}

\rev{Since CAISO is pioneering the design of flexible ramping products, we conduct the same analysis on the WECC 240-bus CAISO model \cite{5299294}. Table \ref{Table:info_39bus} only provides the information for the generators with ramping capacities\footnote{The minimal generation $\underline{g}_n$'s are provided in the model \cite{5299294}.}. Note that, to cope with the increasing renewable energy penetration level, we have tripled the ramping capacities in the system. All the other information are the same as suggested by the model.

The initial energy dispatch is determined by a sample load profile. The total load of the sample profile is 180.15 GW. At time $t=0$ and $t=1$, the forecasted total loads are both 181.15 GW by scaling up the sample load profile. Fig. \ref{fig:WECC_DS} shows the contour of the DS function for this system.

Based on the same goal of renewable energy penetration level, we assume there is a wind plant with a maximal capacity of 80 GW in the system. Fig. \ref{fig:WECC_all} shows the possible distortion costs induced by flexible ramping products for three cases - low wind (16 GW output), modest wind (35 GW output), and high wind (70 GW output). Again, we compare our risk-limiting economic dispatch approach with the greedy approach. As illustrated by Fig. \ref{fig:WECC_all}, our approach performs much better than the greedy approach even when the confidence level is only 80\%. This is because the DS function for this system is even more asymmetric than that for Garver's 6-bus system (see Fig. \ref{fig:contour_6bus}). The asymmetry results in identifying inefficient up and down ramping capacities. As the feasible region shrinks, the greedy approach finally starts getting better. This again highlights that our approach is very promising, since the greedy algorithm does not have any performance guarantee. Numerically, compared with the greedy approach, when the confidence level is greater than 90\%, the average savings are 65.9\% when the predicted wind power is 16 GW, 55.5\% when the predicted wind power is 35 GW, and 56\% when the predicted wind power is 70 GW.}

\subsection{Extensions}
\label{extension}

\begin{figure}[t]
\begin{center}
 \includegraphics[width=7.5cm]{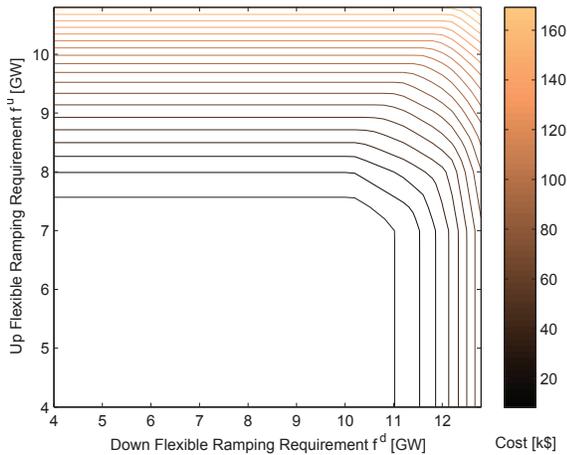}
 %\vspace{-0.3cm}
 \caption{Contour of the DS function for WECC system when $T=3$.}\vspace{-0.2cm}
 \label{fig:WECC_T3}
\end{center}
\end{figure}

\rev{We would like to close this section by briefly discussing the extension of our approach by selecting $T=3$.
To avoid too many arguments, i.e., $F_t^d$'s and $F_t^u$'s, we require the ramping parameters are identical for all time slots, i.e., $F_1^d=F_2^d=F^d$, and $F_1^u=F_2^u=F^u$. Fig. \ref{fig:WECC_T3} shows the contour of the DS function for the WECC system when $T=3$ in this simplified setting.

Comparing Fig. \ref{fig:WECC_DS} and Fig. \ref{fig:WECC_T3}, we may conclude that although the feasible region for $T=3$ (Fig. \ref{fig:WECC_T3}) is smaller than that for $T=2$ (Fig. \ref{fig:WECC_DS}), considering more time slots exploits more constraints in the system - most of the contour lines in Fig. \ref{fig:WECC_T3} have one more segment than the corresponding lines in Fig. \ref{fig:WECC_DS}. Hence, the extension suggests our proposed approach to be even more promising.}

\section{Conclusions and Future Work} \label{conclusions}

This paper proposes a parametric functional analysis of the relations between the generation cost and the key parameters of the flexible ramping products. We present a novel routine to efficiently construct the introduced functions. Such a routine further yields the efficient risk-limiting economic dispatch. Theoretical analysis exhibits valuable information about such relations whereas simulation results further illustrate how such an approach can be used in practice.

This paper can be extended in various directions. For instance, we have not fully investigated the relationship between the distortion cost and the total capacity payment. Also, it is important to analyze the firm behaviors in the electricity markets with the new products: is it easy to gain market power in such a model? \rev{A more accurate empirical evaluation should also include the impacts of existing products, such as frequency regulation, spinning reserves, and non-spinning reserves. Incorporating these existing products will further shrink the feasible region of the MinC function (or equivalently, the DS function), and hence suggest our proposed approach to be even more promising. It is worth noting that our approach is capable of capturing the bidding information for flexible ramping products. Thus, after the bidding structure becomes available in the energy markets, a new empirical evaluation of our approach will be very valuable.}

\rev{In addition, although we have considered the simple generalization to the case when $T=3$ in Section \ref{extension}, the generalization to a more complicated setting, where all the arguments are allowed to be different, is also very interesting and yet more challenging. Odds are that such generalization may exploit more information on how the dynamics of the renewables affect the electricity market.}

\section*{Acknowledgement} \label{ack}

\rev{This research was supported in part by Pennsylvania Infrastructure Technology Alliance and Carnegie Mellon University Scott Institute.}

\bibliographystyle{IEEEtran}
\bibliography{FlexRamp}

%\vspace{-0.2cm}
\appendix

\subsection{Proof of Theorem \ref{thm1}}
\label{app1}

We only prove Part (a) as the proofs for Part (b) and (c) are similar. The continuity and piecewise linearity properties of the DS function are a direct result from \cite[Theorem 1.1-1.3]{Holder2010}. Therefore, we only prove the convexity and monotonicity. Take $f^u$ as an example example and fix $f^d=f^d_0$.

To prove convexity, consider two arbitrary points $$\underline{f}^u(f^d_0)\le f^u_1 < f^u_2 \le \bar{f}^u(f^d_0).$$ Let $g_{n,t}^1$, $r_n^{u,1}$ and $r_n^{d,1}$ denote the optimal solution when solving the optimization problem corresponding to DS($f^u_1,f^d_0$). Similarly, let $g_{n,t}^2$, $r_n^{u,2}$ and $r_n^{d,2}$ denote the optimal solution when solving the optimization problem corresponding to DS($f^u_2,f^d_0$).
For any $f^u(\delta) = \delta f^u_1+(1-\delta)f^u_2,$ where $0\le \delta\le 1$, we can show that
\begin{align}
    g_{n,t}^{\delta} = \delta g_{n,t}^1 + (1-\delta) g_{n,t}^2, \\
    r_n^{u,\delta} = \delta r_n^{u,1} + (1-\delta) r_n^{u,2}, \\
    r_n^{d,\delta} = \delta r_n^{d,1} + (1-\delta) r_n^{d,2},
\end{align}
construct a \emph{feasible} (but not necessarily optimal) solution to the optimization problem corresponding to MinC($f^u(\delta),f^d_0$). Therefore, we have
\begin{equation}\label{chap5:concave}
\begin{aligned}
& \text{MinC}(f^u(\delta),f^d_0) \le  \sum_{t=1}^2 \sum_{n\in\mathcal{N}} C_n g_{n,t}^{\delta}\\
 =&  \delta \sum_{t=1}^2 \sum_{n\in\mathcal{N}} C_n g_{n,t}^1\! +\! (1\!-\!\delta)\sum_{t=1}^2 \sum_{n\in\mathcal{N}} C_n g_{n,t}^2 \\
  =& \delta \ \text{MinC}( f^u_1)+(1-\delta)\ \text{MinC}( f^u_2).
\end{aligned}
\end{equation}
From (\ref{chap5:concave}), for each given $f^d_0$, the corresponding MinC function is a convex function \cite[Section 3.1.1]{boyd2004convex}.

Next, we show that the MinC function is monotonic increasing in the feasible region (\ref{boundfun})-(\ref{boundthetan}) by contradiction.

Suppose the minimum of MinC function is achieved at $f^{u,\star}$ such that $\underline{f}^u(f^d_0)< f^{u,\star} < \bar{f}^u(f^d_0)$. Since the MinC function is continuous, based on the intermediate value theorem, there exists $f^{u,\star} \le \hat{f}^u < \bar{f}^u(f^d_0)$ such that $\text{MinC}(\hat{f}^u,f^d_0)=\text{MinC}(\underline{f}^u(f^d_0),f^d_0)$, which contradicts the definition of $\underline{f}^u(f^d_0)$. From this, and since MinC function is convex, it is monotonic increasing over the range in (\ref{boundfun})-(\ref{boundthetan}).$\hfill \blacksquare$

\vspace{-0.2cm}
\subsection{Proof of Theorem \ref{thm2}}
\label{app2}
Again, we only show the proof for Part (a) since the remaining part is similar.
By the definitions of MinC and MaxUR functions in (\ref{DS_def}) and (\ref{MaxRU}), for any given $f^d_0$, we have
\begin{equation}\label{chap5:MG_NF_1o}
   \text{MinC}(\text{MaxUR}({\theta},f^d_0),f^d_0) \le {\theta},
\end{equation}
\begin{equation}\label{chap5:MG_NF_2o}
\text{MaxUR}(\text{MinC}({f}^u,f^d_0))\ge {f}^u.
\end{equation}
Since the MinC and MaxUR functions are both increasing for given $f^d_0$, their inverses are also increasing. As a result, taking MinC$^{-1}(\cdot)$ at both sides of inequality (\ref{chap5:MG_NF_1o}) leads to
\begin{equation}\label{chap5:MG_NF_1}
    \text{MaxUR}({\theta},f^d_0) \le \text{MinC}^{-1}({\theta},f^d_0).
\end{equation}
Similarly, taking MaxUR$^{-1}$ at both sides of (\ref{chap5:MG_NF_2o}) yields
\begin{equation}\label{chap5:MG_NF_2}
\text{MinC}({f}^u,f^d_0) \ge \text{MaxUR}^{-1}({f}^u,f^d_0).
\end{equation}
By selecting $f^u = \text{MinC}^{-1}({\theta},f^d_0)$, inequality (\ref{chap5:MG_NF_2}) becomes
\begin{equation}\label{chap5:MG_NF_3}
\theta \ge \text{MaxUR}^{-1}(\text{MinC}^{-1}({\theta},f^d_0),f^d_0).
\end{equation}
Using the monotonicity property again, we know
\begin{equation}\label{chap5:MG_NF_3o}
\text{MaxUR}({\theta},f^d_0) \ge \text{MinC}^{-1}({\theta},f^d_0).
\end{equation}
Together, from (\ref{chap5:MG_NF_1}) and (\ref{chap5:MG_NF_3o}), we can conclude that
\begin{equation}\label{chap5:MG_NF_4}
\text{MaxUR}({\theta},f^d_0) = \text{MinC}^{-1}({\theta},f^d_0).
\end{equation}
We can show $\text{MinC}({f}^u,f^d_0) = \text{MaxUR}^{-1}({f}^u,f^d_0)$ in the same way. This concludes the proof.
$\hfill \blacksquare$

\end{document}